\input amstex
\documentstyle{amsppt}
\magnification=\magstep1 \NoRunningHeads
\topmatter
\title
 Finite ergodic index and asymmetry for  infinite measure preserving actions  
 \endtitle

\author
Alexandre I. Danilenko
\endauthor

\email
alexandre.danilenko@gmail.com
\endemail

\address
 Institute for Low Temperature Physics
\& Engineering of National Academy of Sciences of Ukraine, 47 Lenin Ave.,
 Kharkov, 61164, UKRAINE
\endaddress
\email alexandre.danilenko\@gmail.com
\endemail

\abstract
Given $k>0$ and an Abelian countable discrete group $G$ with elements of infinite order, we construct  $(i)$ rigid  funny rank-one infinite measure preserving (i.m.p.) $G$-actions of ergodic index $k$,
$(ii)$ 0-type funny rank-one i.m.p. $G$-actions of ergodic index $k$,
$(iii)$ funny rank-one i.m.p. $G$-actions $T$ of ergodic index 2  such that the product $T\times T^{-1}$  is not ergodic.
It is shown that $T\times T^{-1}$ is conservative for each funny rank-one $G$-action $T$.
\endabstract

 \loadbold

\endtopmatter

\document

\head 0. Introduction
\endhead

Let $G$ be a discrete countable  Abelian group and let $T=(T_g)_{g\in G}$ be a measure preserving action of $G$ on an {\it infinite} $\sigma$-finite standard measure space $(X,\goth B,\mu)$. 
By $T^{-1}$ we denote the ``inverse to $T$'' action of $G$, i.e. $T^{-1}:=(T_{-g})_{g\in G}$. 
Given two $G$-actions $S$ and $R$ of $G$, we denote by $S\times R$ and $S\otimes R$ the following product actions of $G$ and $G\times G$ respectively on the product of the underlying measure spaces:
$S\times R:=(S_g\times R_g)_{g\in G}$ and $S\otimes R:=(S_g\times R_h)_{g,h\in G}$.
If $S$ and $R$ are both conservative or ergodic then $S\otimes R$ is also conservative or ergodic  respectively.
However the product $S\times R$ can be neither ergodic nor conservative.
If $T\times\cdots\times T\,(k\text{ times})$ is ergodic but $T\times\cdots\times T\,(k+1\text{ times})$ is not then $T$ is said to {\it have  ergodic index $k$}.
In 1963, Kakutani and Parry constructed for each $k$,  an infinite Markov shift (i.e. $\Bbb Z$-action) of ergodic index $k$.
In their examples, the product $T\times\cdots\times T\,(k\text{ times})$ is ergodic 
if and only if it is conservative.
For a half a century their examples were the only examples of transformations of finite ergodic index $k>1$.
Recently another family of $\Bbb Z$-actions of an arbitrary finite ergodic index  appeared in \cite{AdSi} by Adams and Silva.
That family  consists of  rank-one transformations $T$ with {\it infinite conservative index}, i.e. $T\times\cdots\times T\,(l\text{ times})$ is conservative for each $l>0$.
We extend and refine their result to the Abelian groups containing  elements of infinite order in the following way.

\proclaim{Theorem 0.1}
Let $G$ has an element of infinite order.
For each $k>0$, there is
a rigid funny rank-one\footnote{Funny rank one  means that there is a sequence $(F_n)_{n=1}^\infty$ of finite subsets in $G$ and a sequence $(A_n)_{n=1}^\infty$ of subsets of finite measure such that $T_gA_n\cap T_gA_n=\emptyset$ whenever $g\ne h\in F_n$ and the sequence of $T$-towers $\{T_gA_n\mid g\in F_n\}$ approximates the entire Borel $\sigma$-algebra as $n\to\infty$.
In the case $G=\Bbb Z^k$, if each $F_n$ is a cube then $T$ is said to be of rank one.} $G$-action  $T$ of ergodic index  $k$.
Moreover,  the $G$-action $\underbrace{T\times\cdots\times T}_{m \text{ times}}\times \underbrace{T^{-1}\times\cdots\times T^{-1}}_{k-m \text{ times}}$ is ergodic for every $m\in\{0,1,\dots,k-1\}$.
Furthermore, if $G=\Bbb Z^d$ for some $d>0$ then $T$ can be chosen in the class of rank-one actions.
\endproclaim

We  note that $T$ has infinite conservative index if $T$ is rigid.
We also note that while the proof of the first claim of Theorem~0.1 in \cite{AdSi} (in the case $G=\Bbb Z$) is somewhat tricky, our proof is shorter  and more direct.

In the next theorem we construct funny rank-one actions of finite ergodic index which are {\it mixing} (called also  {\it  zero type}, see \cite{DaSi} and references therein), i.e. $\lim_{g\to\infty}\mu(T_gA\cap A)=0$ for each subset $A$ of finite measure.
Thus these actions (in the case where $G=\Bbb Z$) are different from those constructed in \cite{KaPa} and \cite{AdSi}.

\proclaim{Theorem 0.2}
Let $G$ has an element of infinite order.
For each $k>0$, there is
a mixing (zero type) funny rank-one $G$-action  $T$ of ergodic index  $k$  such that $T\times\cdots\times T\,(k+1\text{ times})$ is conservative but
$T\times\cdots\times T\,(k+2\text{ times})$ is non-conservative.
Moreover,  the $G$-action $\underbrace{T\times\cdots\times T}_{m \text{ times}}\times \underbrace{T^{-1}\times\cdots\times T^{-1}}_{k-m \text{ times}}$ is ergodic for every $m\in\{0,1,\dots,k-1\}$.
Furthermore, if $G=\Bbb Z^d$ for some $d>0$ then $T$ can be chosen in the class of rank-one actions.
\endproclaim

In a recent paper \cite{Cl--Va}, rank-one transformations $T$ are constructed such that the product $T\times T$ is ergodic but $T\times T^{-1}$ is not.
This is a partial answer to the following question of Bergelson (see problem P10 in \cite{Da1}): is there a transformation $T$ with infinite ergodic index and such that $T\times T^{-1}$ is non-ergodic?
The next theorem extends this result to the actions of Abelian groups and simplifies the original proof.
Moreover, we show (confirming a conjecture from \cite{Cl--Va}) that   these examples do not answer Bergelson's question completely because  the $G$-action $T\times T\times T$ is not even conservative. 

\proclaim{Theorem 0.3}
Let $G$ has an element of infinite order.
 There is a  funny rank-one action  $T$ of $G$ of ergodic index 2 such that  $T\times T^{-1}$ is non-ergodic but conservative  and $T\times T\times T$ is non-conservative.
Furthermore, if $G=\Bbb Z^d$ for some $d>0$ then $T$ can be chosen in the class of rank-one actions.
\endproclaim

It follows, in particular, that $T$ is asymmetric, i.e. not isomorphic to $T^{-1}$. 
Thus, Theorem~0.3 illustrates that even such a simple invariant as ``ergodicity of products''  can distinguish between $T$ and $T^{-1}$.
Of course, this is impossible in the framework of finite measure preserving actions.
For other, more involved  examples of asymmetric infinite measure preserving systems   we refer to
\cite{DaRy} and \cite{Ry}.

\comment
We also prove a counterpart  of Theorem~0.3.

\proclaim{Theorem 0.4}
Let $G$ has an element of infinite order.
 There is a  funny rank-one action  $T$ of $G$ such that $T\times T$ is non-ergodic but conservative, $T\times T\times T$  is non-conservative and $T\times T^{-1}$ is ergodic.
Furthermore, if $G=\Bbb Z^d$ for some $d>0$ then $T$ can be chosen in the class of rank-one actions.
\endproclaim

\endcomment

It was shown in \cite{Cl--Va} that for each rank-one $\Bbb Z$-action $T$, the product $T\times T^{-1}$ is conservative.
We generalize this result to the funny rank-one action of Abelian groups.

\proclaim{Theorem 0.4}
Let $T$ be a funny rank-one action of $G$.
Then the $G$-action $T\times T^{-1}$ is conservative.
\endproclaim

On the other hand, we show that  each infinite measure preserving Markov shift $T$ of ergodic index 1,
the product $T\times T^{-1}$ is not conservative (Corollary~3.3).
This was also proved in \cite{Cl-Va} under an additional assumption that $T$ is ``reversible'' as a Markov shift.
It follows from Corollary~3.2 that if an infinite Markov shift $T$ has an ergodic index higher than 1 then $T\times T^{-1}$ is ergodic.
Hence within the class of infinite Markov shifts, the answer to Bergelson's question is negative.

\head 1. $(C,F)$-construction
\endhead

All the examples in this paper are built via the $(C,F)$-construction which is an algebraic counterpart of the classical ``geometric'' cutting-and-stacking inductive construction process with a single tower on each step.
In this section we  briefly  outline the basics of this construction.
For a detailed exposition we refer the reader to \cite{Da1} and \cite{Da2}.
Given two finite subsets $A,B\subset G$, we denote by $A+B$ the set of all sums $\{a+b\mid a\in A,b\in B\}$.  If $A$ is a singleton, say $A=\{a\}$, we write $a+B$ in place of $\{a\}+B$. 

Let $(F_n)_{n\ge 0}$  and $(C_n)_{n\ge 1}$ be two sequences of finite subsets in $G$ such that for each $n>0$,
$$
\align
&F_0=\{0\},  \quad \# C_n>1,\tag 1-1\\
&F_n+C_{n+1}\subset F_{n+1},\tag1-2\\
&(F_n+c)\cap (F_n+c')=\emptyset\quad\text{if $c,c'\in C_{n+1}$ and $c\ne c'$.}\tag1-3
\endalign
$$
We let $X_n:=F_n\times C_{n+1}\times C_{n+2}\times\cdots$ and endow this set with the infinite product topology.
Then $X_n$ is a compact Cantor (i.e. totally disconnected  perfect metric) space.
The  mapping
$$
X_n\ni (f_n,c_{n+1},c_{n+2},\dots)\mapsto(f_nc_{n+1}, c_{n+2},\dots)\in X_{n+1}
$$
is a topological embedding
of $X_n$ into $X_{n+1}$.
Therefore an inductive limit $X$ of the sequence $(X_n)_{n\ge 0}$ furnished with these embeddings  is a well defined locally compact Cantor  space.
Given a subset $A\subset F_n$, we let 
$$
[A]_n:=\{x=(f_n,c_{n+1},\dots)\in X_n\mid f_n\in A\}
$$ 
and call this set an  {\it $n$-cylinder} in $X$.
It is open and compact in $X$.
The collection of all cylinders coincides with the family of  all compact open subset  in $X$.
It is easy to see that
$$
\align
[A]_n\cap[B]_n &=[A\cap B]_n, \quad
[A]_n\cup[B]_n =[A\cup B]_n\quad \text{and}\\
 [A]_n &=[A+C_{n+1}]_{n+1}
 \endalign
 $$
  for all $A,B\subset F_n$ and  $n\ge 0$.
For brevity, we will write $[f]_n$ for $[\{f\}]_n$, $f\in F_n$.
Now we define the
 {\it $(C,F)$-measure} $\mu$ on $X$ by setting
$$
\mu([A]_n)=\frac{\# A}{\# C_1\cdots\# C_n} \quad\text{for each subset }A\subset F_n, n>0,
$$
and then extending $\mu$ to the Borel $\sigma$-algebra  on $X$.
We note that $\mu$ is infinite if and only if
$$
\lim_{n\to\infty}\frac{\# F_n}{\# C_1\cdots\#C_n}=\infty.
\tag1-4
$$
Suppose that for each $g\in G$,
 $$
g+F_n+C_{n+1}\subset F_{n+1}\quad\text{eventually in $n$.}\tag1-5
 $$
We now define an action of $G$ on $X$.
Given $x\in X$ and $g\in G$, there is $n>0$ such that  $x=(f_n,c_{n+1},c_{n+2},\dots)\in X_n$ and $g+f_n\in F_n$.
We  let 
$$
T_gx:=(g+f_n, c_{n+1},\dots)\in X_n\subset X.
$$
Then $T_g$ is a well defined homeomorphism of $X$ and 
 $T:=(T_g)_{g\in G}$ is a continuous action of $G$ on $X$.
We call it {\it the $(C,F)$-action of $G$ associated with the sequence $(C_n,F_{n-1})_{n\ge 1}$}.
It is free.
If $x,y\in X_n$, $x=(f_n,c_{n+1,\dots})$, $x'=(f',c_{n+1}',\dots)$ and $y=T_gx$ then $g=(f'-f)+(c_{n+1}'-c_{n+1})+\cdots$.
Only  finitely many parenthesis in this infinite sum are different from $0$.
We note that $T$ is of {\it funny rank one along $(F_n)_{n\ge 0}$}, because the sequence of $F_n$-towers $\{T_f[0]_n\mid f\in F_n\}=\{[f]_n\mid f\in F_n\}$ approximates the entire Borel $\sigma$-algebra on $X$ as $n\to\infty$.
It is easy to see that  $T$ preserves  $\mu$.
We note that   the action $T\otimes T$ of $G\times G$ is   also a $(C,F)$-action.
It is associated with the sequence $(C_n\times C_n,F_{n-1}\times F_{n-1})_{n\ge 1}$.
Therefore if $A$ is a subset of $F_n\times F_n$ then we have $[A]_n=\bigsqcup_{(a,b)\in A}[a]_n\times[b]_n$.

To state the following lemma we recall the definition of  {\it full grouppoid}.
Given an measure preserving action $R=(R_g)_{g\in G}$ on a standard measure space $(X,\mu)$ and a subset $A\subset X$,
we say that  a Borel map $\tau:A\to X$ belongs to the full grouppoid of $R$ (and write $\tau\in[[R]]$) if $\tau$
is one-to-one and $\tau x\in\{R_gx\mid g\in G\}$ for all $x\in A$.
Equivalently, there is a partition of $A$ into subsets $A_g$, $g\in G$,  such that $\tau x=R_gx$ if $x\in A_g$
and $R_gA_g\cap R_hA_h=\emptyset$ if $g\ne h\in G$.
Some of $A_g$ can be of zero measure. 
It follows that $\tau$ preserves $\mu$.

\proclaim{Lemma 1.1}
 Let $\delta>0$,  let $H$ be a subgroup of $G$ and let $\Cal N$ be an infinite subset of $\Bbb N$.
 \roster
 \item"\rom(i)"
 If for each $n\in\Cal N$, there is a subset $A\subset[0]_n$ and a  map $\tau: A\to[0]_n$ such that  $\tau\in[[(T_h)_{h\in H}]]$,
  $\mu(A)\ge \delta\mu([0]_n)$  and $\tau x\ne x$ for all $x\in A$
then the action $(T_h)_{h\in H}$ is conservative.
  \item"\rom(ii)"
If for each $n\in\Cal N$ and $v,w\in F_n$, there is a subset $A\subset[v]_n$ and a map $\tau: A\to[w]_n$ 
such that $\tau\in[[(T_h)_{h\in H}]]$
and $\mu(A)\ge \delta\mu([v]_n)$  then the action $(T_h)_{h\in H}$ is ergodic.
\endroster
\endproclaim

\demo{Proof}
(i) Let $B$ be a subset of $X$ of positive measure.
Then there is $n\in\Cal N$ and $f\in F_n$ with $\mu([f]_n\cap B)>(1-\frac\delta 4)\mu([f]_n)$.
By the assumption of the lemma, there is a subset $A\subset[0]_n$ and  a  map $\tau:A\to[0]_n$ such that $\tau\in[[(T_h)_{h\in H}]]$, $\mu(A)>\delta\mu([0]_n)$ and $\tau x\ne x$ for all $x\in A$.
We define a new map $\varphi:T_fA\to[f]_n$ by setting $\varphi:=T_f\tau T_f^{-1}$.
Since $G$ is Abelian, $\varphi\in[[(T_h)_{h\in H}]]$.
Moreover, $\varphi x\ne x$ for all $x\in T_fA$ and
$$
\mu(\varphi(T_fA\cap B)\cap B)>\frac\delta 2\mu([f]_n)>0.
$$
Therefore, there is $h\in H$ such that $h\ne 0$ and $\mu(T_h(T_fA\cap B)\cap B)>0$.
Hence $(T_h)_{h\in H}$ is conservative.

(ii) Let $B_1$ and $B_2$ be  subsets of $X$ of positive measure.
Then there are $n\in\Cal N$ and elements $v,w\in F_n$ with $\mu([v]_n\cap B_1)>(1-\frac\delta 4)\mu([v]_n)$ and $\mu([w]_n\cap B_2)>(1-\frac\delta 4)\mu([w]_n)$.
By the assumption of the lemma, there is a subset $A\subset[v]_n$ and   a  map $\tau:A\to[w]_n$ such that $\tau\in[[(T_h)_{h\in H}]]$ and  $\mu(A)>\delta\mu([v]_n)$. 
It follows that $\mu(\tau(B_1\cap[v]_n)\cap[w]_n\cap B_2)>0$.
Therefore there is $h\in H$ such that $\mu(T_hB_1\cap B_2)>0$.
Hence $(T_h)_{h\in H}$ is ergodic.
\qed
\enddemo


\head 2. Proof of the main results
\endhead

Fix a F{\o}lner sequence $(\Cal F_n)_{n=1}^\infty$ in $G$ such that $\Cal F_1\subset \Cal F_2\subset\cdots$  and $\bigcup_n\Cal F_n=G$.
In the case where $G=\Bbb Z^d$, we choose $\Cal F_n$ to be a cube for each $n$.
The actions whose existence is stated in Theorems~0.1--03 will appear as  $(C,F)$-actions associated with some sequences $(C_n,F_{n-1})_{n\ge 1}$.
Moreover, $(F_n)_{n=1}^\infty$ will be a subsequence of $(\Cal F_n)_{n=1}^\infty$.
Therefore in the case $G=\Bbb Z^d$, the associated $(C,F)$-actions will be automatically of rank one. Hence we do not need to prove the final claims of Theorems~0.1--0.3.

 \demo{Proof of Theorem~0.1} (i)
 Partition the natural numbers $\Bbb N$ into countably many  subsets $\Cal N_f$ indexed by elements $f\in G^k$ such that every  $\Cal N_f$ is an infinite arithmetic sequence. 
 For each $f=(f_1,\dots,f_k)\in G^k$ and each $n\in\Cal N_f$ there is a unique sequence $(d_{n,j})_{j=0}^k$,  of nonnegative integers such that $d_{n,0}=0$ and $d_{n,j-1}-d_{n,j}=f_{j}$ for all $j=1,\dots,k$.
 Fix an increasing sequence $(R_n)_{n\ge 0}$ of positive integers such that
$\sum_{n\ge 0}R_n^{-k}=+\infty$ but $\sum_{n\ge 0}R_n^{-k-1}<+\infty$.
We note that then $\sum_{n\in\Cal N_f}R_n^{-k}=+\infty$ for each $f\in G^k$.

 To construct $T$ we have to define the corresponding sequence $(C_n,F_{n-1})_{n\ge 1}$.
 This will be done inductively.
 We let $F_0=\{0\}$.
Suppose that we have already determined  the sequence $(C_j,F_{j})_{j=1}^{n-1}$.
\comment
 Select a sequence $(d_{n,1},\dots, d_{n,k})_{n=1}^\infty$ in $G^k$ such that each element of $G^k$ occurs infinitely often in this sequence.
 \endcomment
 Then we let 
 $$
 C_{n,0}:=\{0,a_{n}+d_{n,1},\dots,ka_{n}+d_{n,k}\},\quad C_{n,1}:=\{w_{n},2w_{n},\dots,(R_n-k-1)w_n\},
 $$
 and $C_n:=C_{n,0}\sqcup C_{n,1}$, where
 the elements $a_{n},w_{n}\in G$,  are chosen so that
 $$
 (C_{n,1}-C_n-C_n+C_n)\cap(F_{n-1}-F_{n-1}+F_{n-1}-F_{n-1})=\{0\}.\tag2-1
 $$
  Now let
$F_n$ be an element of  $(\Cal F_j)_{j\ge 1}$ such that 
$
F_{n-1}+F_{n-1}+C_n\subset F_n.
$
Continuing this  process infinitely many  times we obtain a sequence $(C_n,F_{n-1})_{n\ge 1}$ satisfying~\thetag{1-1}--\thetag{1-5}.
Denote by $T$ the associated $(C,F)$-action of $G$.
Let $(X,\mu)$ stand for the space of $T$.
It is easy to see that $\mu(T^{w_n}A\cap A)\to \mu(A)$ as $n\to\infty$ for each subset $A$ of finite measure.
Hence $T$ is rigid.
 
 {\bf Claim 1.}
 $T\times\cdots\times T\,(k\text{ times})$ is ergodic.
 
 Take $n>0$ and let  $v,w\in F^k_n$.
 We let $f:=w-v$.
 Given  $x\in[v]_n\subset X^k$, we write  the expansion
 $$
 x=(v,x_{n+1},x_{n+2},\dots)\in F_n^k\times C_{n+1}^k\times C_{n+2}^k\times\cdots
 $$ 
 and set
 $$
 \ell(x):=\min\{l\in\Cal N_f\cap\{n+1,n+2,\dots\}\mid 
 x_l\in (C_{l,0}\setminus\{ka_l+d_{l,k}\})^k\}.
 $$
 Let $A$ denote  the subset of $[v]_n$ where the map $\ell$ is well defined.
 Then
 $$
\frac{\mu^k([v]_n\setminus A)}{\mu^k([v]_n)}= \prod_{l\in\Cal N_f,l>n}\frac{(\# C_l)^ k-(\# C_{l,0}-1)^k} {(\# C_l)^k}=
  \prod_{l\in\Cal N_f,l>n}\bigg(1-\frac{k^k}{R_l^k}\bigg)=0
 $$
 because $\sum_{l\in\Cal N_f}{R_l^{-k}}=\infty$. 
 Thus $\ell$ is defined almost everywhere on  $[v]_n$.
 For $l>n$, we set
 $$
 A_l:=\{x\in A\mid \ell(x)=l\text{ and } x_l=(0,a_l+d_{l,1},\dots,(k-1)a_l+d_{l,k-1})\}.
 $$
 Then $\mu^k(\bigsqcup_{l>n}A_l)=\mu^k(A)/(k+1)^k$.
 We now define a map
 $\tau: \bigsqcup_{l>n}A_l\to X^k$ by setting
 $$
 \tau x:=(T_{a_l}\times\cdots\times T_{a_l})x\quad\text{if }x\in A_l, \, l>n.
 $$
Of course, $\tau\in[[ T\times\cdots\times T ]]$.
 Since 
 $$
 a_l+ ((j-1)a_l+d_{l,j-1})=f_{j}+(ja_l+d_{l,j})\quad \text{for $j=1,\dots,k$},
 $$
  it follows that 
 $(T_{a_l}\times\cdots\times T_{a_l})x=(T_{f_1}\times\cdots\times T_{f_k})y$, where
 $y=(y_i)_{i\ge n}\in F_n^k\times C_{n+1}^k\times C_{n+2}^k\times\cdots$,  $y_i=x_i$ if $i\ge n$ and $i\ne l$ and 
 $y_l=(a_l+d_{l,1},\dots,ka_l+d_{l,k})\in C_l^k$.
  Since $y\in[v]_n$ and $f=w-v$, we obtain that $(T_{f_1}\times\cdots\times T_{f_k})y\in[w]_n$ for each $x\in A_l$.
  Hence 
 $\tau x\in [w]_n$ for each $x\in \bigsqcup_{l>n}A_l$.
Therefore  $T\times\cdots\times T\,(k\text{ times})$ is ergodic by Lemma~1.1(ii).

{\bf Claim 2.}
The $G$-action $\underbrace{T\times\cdots\times T}_{m \text{ times}}\times \underbrace{T^{-1}\times\cdots\times T^{-1}}_{k-m \text{ times}}$ is ergodic for every $m\in\{0,1,\dots,k-1\}$.
The proof is  similar to the proof of Claim 1.
There are only a few points of difference which we specify now. 
Let $\widetilde f:=(f_1,\dots,f_m,-f_{m+1},\dots,-f_k)$. 
Replace $\Cal N_f$ with $\Cal N_{\widetilde f}$ in the definition of $\ell$.
Define
$$
\multline
 B_l:=\{x\in A\mid \ell(x)=l\text{ and } 
 x_l=(0,a_l+d_{l,1},\dots,(m-1)a_l+d_{l,m-1},\\(m+1)a_l+d_{l,m+1},\dots,
 ka_l+d_{l,k})\}.
 \endmultline
 $$
Replace $A_l$ with $B_l$ and $T_{a_l}\times\cdots\times T_{a_l}$
with $\underbrace{T_{a_l}\times\cdots\times T_{a_l}}_{m \text{ times}}\times T_{-a_l}\times\cdots\times T_{-a_l}$ in the definition of $\tau$.

 {\bf Claim 3.}
 $T\times\cdots\times T\,(k+1\text{ times})$ is not ergodic.
 Choose $n>0$ such that $\sum_{j=n}^\infty(\frac{k+1}{R_j})^{k+1}<1$.
 We now let
 $$
 W:=\{\bold 0\}\times (C_n^{k+1}\setminus C_{n,0}^{k+1})\times (C_{n+1}^{k+1}\setminus C_{n+1,0}^{k+1})\times\cdots
\subset[\bold 0]_{n-1}\subset X^{k+1}.
 $$
Here $\bold 0$ denotes zero in  $G^{k+1}$.
 Then 
 $$
\frac{ \mu^{k+1}(W)}{\mu^{k+1}([\bold 0]_{n-1})}=
\prod_{j\ge n}\bigg(1-\bigg(\frac{k+1}{R_j}\bigg)^{k+1}\bigg)\ge
1-\sum_{j\ge n}\bigg(\frac{k+1}{R_j}\bigg)^{k+1}>0.
 $$
Fix $h\in F_{n-1}\setminus\{ 0\}$.
 We now show that the $(T\times\cdots\times T)$-orbit of $W$ does not intersect the cylinder $B:= [0]_{n-1}\times \cdots\times [0]_{n-1}\times[h]_{n-1}\subset X^{k+1}$.
 If not, then there is $x=(x^1,\dots,x^{k+1})\in W$ and $g\in G$ such that 
 $$
 (T_gx_1,\dots,T_gx_{k+1})\in B.\tag2-2
 $$
 Consider the expansions 
 $$
 \align
 x^l&=(0,x^l_{n},x^l_{n+1},\dots)\in \{0\}\times C_n\times C_{n+1}\times\cdots,\quad l=1,\dots,k+1,\\
T_gx^l&=(0,y^l_{n},y^l_{n+1},\dots)\in \{0\}\times C_n\times C_{n+1}\times\cdots,\quad l=1,\dots,k,\quad\text{and}\\ 
T_gx^{k+1}&=(h,y^{k+1}_{n},y^{k+1}_{n+1},\dots)\in \{h\}\times C_n\times C_{n+1}\times\cdots.\\ 
 \endalign
 $$
It follows from \thetag{2-2} that there are integers $N_l\ge n$, $l=1,\dots,k+1$, such that
 $$
\cases
 g=\sum_{i=n}^{N_l}(y_i^l-x_i^l),\quad l=1,\dots,k,\\
  g=h+\sum_{i=n}^{N_{k+1}}(y_i^{k+1}-x_i^{k+1})
 \endcases
 \tag2-3
 $$
 and
 $y_{N_l}^l\ne x_{N_l}^l$, $l=1,\dots,k+1$.
Then  \thetag{2-1} yields that $N_1=\cdots=N_{k+1}$.
 Since $x\in W$, there exists $l_0\in\{1,\dots,k+1\}$ with $x^{l_0}_{N_1}\in C_{N_1, 1}$.
 It now follows from \thetag{2-1} that $y_{N_1}^l-x_{N_1}^l=y_{N_1}^{l_0}-x_{N_1}^{l_0}$ for all $l=1,\dots,k+1$.
 Hence we deduce from \thetag{2-3} that
 $$
 \cases
 g-(y_{N_1}^1-x_{N_1}^1)=\sum_{i=n}^{N_1-1}(y_i^l-x_i^l),\quad l=1,\dots,k,\quad\text{and}\\
  g-(y_{N_1}^1-x_{N_1}^1)=h+\sum_{i=n}^{N_{1}-1}(y_i^{k+1}-x_i^{k+1}).
 \endcases
 $$
Repeating this procedure  at most $N_1-n-1$ times we obtain that $g=g-h$, a contradiction.
\qed
 \comment
Given $c,c'\in C_n$, we write $c\not\approx c'$ if $c=b_{n,i}+a_{n,j}+d_{n,j}$
 and $c'=b_{n,i'}+a_{n,j'}+d_{n,j'}$ with $i\not =i'$.
 Let 
 $$
 W:=\{(x^1,\dots,x^{k+1})\in [0]_1^{k+1}\mid x^i_l\not\approx x^j_l\text{ whenever }i\ne j\text{ for each }l>1\},
 $$
 where $x^j_l\in C_l$ is the $l$-th coordinate of $x^j\in[0]_1$.
 Then 
 $$
 \mu^{k+1}(W)=\frac{1}{(\# C_1)^{k+1}}\prod_{n=2}^\infty\bigg(1-\frac 1{N_n}\bigg)\cdots\bigg(1-\frac{k}{N_n}\bigg)>0.
 $$
 Fix $h\in F_1\setminus\{0\}$.
  We claim that the   $T\times\cdots\times T\,(k+1\text{ times})$-orbit of $W$ does not intersect the cylinder $[0]_1\times\cdots\times [0]_1\times[h]_1\subset X^{k+1}$.
 If not, there is a point $(x^1,\dots,x^{k+1})\in W$ and an element $g\in G$ such that 
 $$
 \aligned
 &T_gx^j\in [0]_1,\quad j=1,\dots,k,\quad\text{and}\\
  &T_gx^{k+1}\in [h]_1.
 \endaligned
 \tag2-9
$$
Denote by $y^j_l\in C_l$ the $l$-th coordinate of 
 $T_gx^j$, $l>1$.
Then \thetag{2-9} yields that there are integers $M_1,\dots,M_{k+1}$ such that
$$
\aligned
g&=(y^j_2-x^j_2)+\cdots+ (y^j_{M_i}-x^j_{M_i}),\quad j=1,\dots,k,\quad\text{and}\\
g-h&=(y^{k+1}_2-x^{k+1}_2)+\cdots+ (y^{k+1}_{M_{k+1}}-x^{k+1}_{M_{k+1}}),
\endaligned
\tag2-10
$$
with $y^j_{M_i}-x^j_{M_i}\ne 0$ for each $j=1,\dots,k+1$.
From this and \thetag{2-1} we deduce that $M_1=M_2=\cdots =M_{k+1}$.
We now write 
$$
x^j_{M_1}=b_{M_1,i_j}+a_{M_1,l_j}+d_{M_1,l_j}\quad\text{and}\quad 
y^j_{M_1}=b_{M_1,s_j}+a_{M_1,r_j}+d_{M_1,r_j}
$$ 
for some $i_j,s_j\in\{1,\dots,N_n\}$ and $l_j,r_j\in\{1,\dots,k\}$, $j=1,\dots,k+1$.
It now follows from  \thetag{2-10} that 
$$
\align
\{s_j,i_1\}&=\{i_j,s_1\}\quad\text{and}\tag2-11\\
 \{r_j,l_1\}&=\{l_j,r_1\}\tag2-12
\endalign
$$ 
for each $j=2,\dots,k+1$.
Since $i_1\ne i_j$ for each $j\ne 1$ by the definition of $W$, \thetag{2-11} yields that $i_j=s_j$ for each $j=1,\dots,k+1$.
It follows from \thetag{2-12} that for each $j$, we have either $r_j=l_j$ or $r_j=r_1$ and $l_j=l_1$.
If the first condition happens then $x^j_{M_1}=y^j_{M_1}$, that contradicts to the choice of  $M_1$.
Therefore $r_j=r_1$ and $l_j=l_1$ for each $j$.
Thus we deduce from \thetag{2-10} that
$$
\aligned
g'
&=(y^j_2-x^j_2)+\cdots+ (y^j_{M_1-1}-x^j_{M_1-1}),\quad j=1,\dots,k,\quad\text{and}\\
g'-h&=(y^{k+1}_2-x^{k+1}_2)+\cdots+ (y^{k+1}_{M_{1}-1}-x^{k+1}_{M_{1}-1}),
\endaligned
$$
where $g':=g+a_{M_1,r_1}-d_{M_1,r_1}-a_{M_1,l_1}+d_{M_1,l_1}$.
Repeating this procedure at most $M_1-2$ times we obtain that $g=g-h$, a contradiction.

 (ii) 
 We  set
 $C_n':=C_n\sqcup\{B_n,2B_n,\dots,N_nB_n\}$, where $B_n$ is a huge number.
 Let $T$ be the $(C,F)$-action associated with $(C_n',F_{n-1})_{n\ge 1}$.
 The proof that $T\times\cdots\times T(k\text{ times})$ is ergodic is similar to the proof of Claim~1.
 Moreover, it is easy to see that $T$ is partially rigid, namely $\lim_{n\to\infty} \mu(T_{B_n}A\cap A)\ge{\mu(A)}/{(k+2)}$ for each subset $A$ of positive measure.
 Hence all Cartesian powers of $T$ are conservative.
 It remains to show that the $(k+1)$-power of $T$ is non-ergodic.
 
 Otherwise there is $g\in G$ and $M>0$ such that
 $$
 \align
g&=(y^i_1-x^i_1)+\cdots+ (y^i_{M}-x^i_{M}),\quad i=1,\dots,k,\\
g+h&=
(y^{k+1}_1-x^{k+1}_1)+\cdots+ (y^{k+1}_{M}-x^{k+1}_{M})
\endalign
$$
and
$y^i_{M}\ne x^i_{M}$ for each $i=1,\dots,k+1$.
It follows that  elements $y^i_M-x^i_M\in (C_M'-C_M')\setminus\{0\}$, $i=1,\dots,k+1$, are close to each other.
Then one of the following three cases happens:

(a) $y^i_M,x^i_M\in C_M$ for all $i=1,\dots,k+1$.

(b) $y^i_M,x^i_M\in\{B_M,2B_M,\dots,N_MB_M\}$ for all $i=1,\dots,k+1$.

(c) there is $i_0$ such that $y_{i_0}\in C_M$ and  $x_{i_0}\in \{B_M,2B_M,\dots,N_MB_M\}$.

It was shown in the proof of (i) that the case (a) is impossible.
In the case (b), we have $y^1_M-x^1_M=\cdots =y^{k+1}_M-x^{k+1}_M$.
Hence we can ``reduce'' $M$ by 1 by considering $g-y^1_M+x^1_M$ instead of $g$.
The same we have in the case (c).
Thus repeating  $M$ times, we obtain that $g=g+h$, a contradiction. 
\qed

\endcomment

 \enddemo

\demo{Proof of Theorem~{0.2}}
The desired action is constructed in the same way as in Theorem~1 however $C_{n,1}$ is different.
We now set
$$
C_{n,1}:=\{w_{n,1},\dots,w_{n,R_n-k-1}\},
$$
where the elements $w_{n,j}\in G$ are chosen so that \thetag{2-1} is satisfied and
$$
\text{the mapping \ $(C_{n,1}\times C_n)\setminus\Cal D\ni (c,c')\mapsto c-c'\in G$ is one-to-one},\tag2-4
$$
where $\Cal D$ is the diagonal in $G\times G$.
As in the proof of Theorem~0.1, we denoted the corresponding $(C,F)$-action by $T$.
Claims 1--3 from the proof of that theorem hold (verbally) for the ``new'' $T$ as well.

{\bf Claim 4.} $T\times\cdots\times T\,(k+1\text{ times})$ is conservative.

Take $n>0$.
Given $x=(x^1,\dots,x^{k+1})\in [\bold 0]_{n}$, 
we set
$$
\ell(x):=\min\{l> n\mid x_l^1=\cdots=x_l^{k+1}\}.
$$
Let $A$ denote the subset of $[\bold 0]_{n}$ where $\ell$ is well defined.
Then
$$
\frac{\mu^{k+1}([0]_n\setminus A)}{\mu^{k+1}([0]_n)}= \prod_{l>n}\frac{(\# C_l)^ {k+1}-\# C_{l}} {(\# C_l)^{k+1}}=
  \prod_{l\in\Cal N_f,l>n}\bigg(1-\frac 1{R_l^k}\bigg)=0.
 $$
 For each $m\in\Bbb N$ and $c\in C_m$, we let $A_{m,c}:=\{x\in A\mid \ell(x)=m, x_m^1=c\}$ and fix an element $c'$
 from $C_m\setminus\{c\}$.
We now define a map $\tau:A\to X^{k+1}$ by setting
$$
\tau x=(T_{c'-c}\times\cdots\times T_{c'-c} )x\quad\text{if }x\in A_{m,c},\ c\in C_m,\ m\in\Bbb N.
$$
Since $A=\bigsqcup_{m\in\Bbb N}\bigsqcup_{c\in C_m}A_{m,c}$, it follows that
$\tau x$ is well defined,  $\tau x\in [\bold 0]_n$ and 
 $\tau\in[[T\times\cdots\times T]]$.
It remains to apply Lemma~1.1(i). 

{\bf Claim 5.} $T\times\cdots\times T\,(k+2\text{ times})$ is not conservative.

Choose $n>0$ such that $\sum_{j=n}^\infty  R_j^{-k-1}<0.5$ and $R_n>(k+1)^{k+2}$.
Denote by $D_n$ the diagonal in $C_{n,1}^{k+2}$, i.e. $D_n:=
\{(c_1,\dots,c_{k+2})\in C_{n,1}^{k+2}\mid c_1=\cdots=c_{k+2}\}$.
 We now let
 $$
 W:=(C_n^{k+2}\setminus (C_{n,0}^{k+2}\cup D_{n}))\times (C_{n+1}^{k+2}\setminus (C_{n+1,0}^{k+2}\cup D_{n+1}))\times\cdots
\subset[\bold 0]_{n-1},
 $$
 where $\bold 0$ stands now for the zero in $G^{k+2}$.
 Then 
 $$
 \align
\frac{ \mu^{k+2}(W)}{\mu^{k+2}([\bold 0]_{n-1})}&=
\prod_{j\ge n}\bigg(1-\bigg(\frac{k+1}{R_j}\bigg)^{k+2}-\frac{R_j-k-1}{R_j^{k+2}}\bigg)\\
&\ge
\bigg(1-\sum_{j\ge n} \frac 1{R_j^{k+1}}\bigg(1+\frac{(k+1)^{k+2}}{R_j}\bigg)\bigg)>0.
\endalign
 $$
 We now show that  $W$ 
 is a $(T\times\cdots\times T)$-wandering set.
 If not, then there is $x=(x^1,\dots,x^{k+2})\in W$ and $g\in G$ such that 
 $
 (T_gx^1,\dots,T_gx^{k+2})\in W.
 $
 Consider the expansions 
 $$
 \align
 x^l&=(0,x^l_{n},x^l_{n+1},\dots)\in \{0\}\times C_n\times C_{n+1}\times\cdots\quad \text{and}\\
T_gx^l&=(0,y^l_{n},y^l_{n+1},\dots)\in \{0\}\times C_n\times C_{n+1}\times\cdots,
 \endalign
 $$
 $l=1,\dots,k+2$.
Arguing in the same way as in the proof of Claim~3, we find $N_1$ such that
 $
 g=\sum_{i=n}^{N_1}(y_i^l-x_i^l),
 $
  $0\ne y_{N_1}^l-x_{N_1}^l=y_{N_1}^{1}-x_{N_1}^{1}$   for all $l=1,\dots,k+2$.
  Moreover,
  $x^{l}_{N_1}\in C_{N_1,1}$ for all $l=1,\dots,k+2$ because  $x\in W$.
Then we deduce from~\thetag{2-4} that $x_{N_1}^1=\cdots=x_{N_1}^{k+2}$, i.e. $(x_{N_1}^1,\dots,x_{N_1}^{k+2})\in D_{N_1}$.
 Therefore $x\notin W$, a contradiction.

{\bf Claim 6.} $T$ is mixing.
Let $A\subset F_{n-1}$, $g\in (F_{n}-F_n)\setminus (F_{n-1}-F_{n-1})$ and $g+ A+C_n\subset F_n$.
Then we have
$$
\align
\mu(T_g[A]_{n-1}\cap [A]_{n-1})
&
=\sum_{c,c'\in C_{n}}\mu([g+A+c]_{n}\cap[A+c']_n)\\
&\le\frac{\mu([A]_{n-1})\#\{(c,c')\in C_n\times C_n\mid g\in A-A+c-c'\}}{\# C_n}.
\endalign
$$
We first note that if 
$$
g\in A-A+c-c'\tag2-5
$$ 
then $c\ne c'$ by the choice of $g$. 
If either $c\in C_{n,1}$ or $c'\in C_{n,1}$ then we deduce from~\thetag{2-1} and \thetag{2-4}  that at most one such pair $(c,c')$ satisfies \thetag{2-5}.
  If both $c\not\in C_{n,1}$ and $c'\not\in C_{n,1}$ then $c,c'\in C_{n,0}$ and hence there are no more than 
  $k(k+1)$ such pairs $(c,c')$ satisfying \thetag{2-5}.
  Hence 
$$
\mu(T_g[A]_{n-1}\cap [A]_{n-1})< \frac{(k+1)^2}{\# C_n}\mu( [A]_{n-1}).
$$
It follows that $\lim_{g\to\infty}\mu(T_gB\cap B)=0$ for each cylinder $B$ and hence for each subset of finite measure in $X$.
\qed
\enddemo

\demo{Proof of Theorem~{0.3}}
Let  $(d_n)_{n=1}^\infty$ be a sequence of elements of $G$ in which each element of $G$ occurs infinitely many times.
Let $(N_n)_{n=1}^\infty$ be a sequence of positive integers such that
  $\sum_{n>0}\frac 1{N_n}<\frac 14$.

Suppose that we have already determined   $(C_j,F_{j})_{j=1}^{n-1}$.
Suppose also that $d_n\in F_{n-1}-F_{n-1} $.
We then  set
$$
C_n:=\{e_{n,i},-e_{n,i},-l_{n,i}, l_{n,i}-d_n\mid i=1,\dots N_n\},
$$
for some elements $e_{n,i},  l_{n,i}$ of $G$, $1\le i
\le N_n$ such that
$$
(C_n-C_n)\cap(F_{n-1}-F_{n-1}+F_{n-1}-F_{n-1})=\{0\}.\tag2-6
$$
We call
 $e_{n,i}$ and $-e_{n,i}$ as well as  $l_{n,i}$ and $-l_{n,i}-d_n$ {\it antipodal}, $1\le i\le N_n$.
 If $c_1,\dots,c_4\in C_n$,   $c_1$ and $c_4$ are antipodal and $c_2$ and $c_3$ are antipodal then  
 $$
 (c_1-c_2)-(c_3-c_4)\in \{0,d_n,-d_n\}.
 $$
We introduce the following conditions on $C_n$.
Let $c_1,\dots,c_4\in C_n$.
$$
\gather
{\aligned
&\text{If  $0\ne (c_1-c_2)-(c_3-c_4)\in F_{n-1}-F_{n-1}+F_{n-1}-F_{n-1}$}\\
&\text{ then $c_1$ and $c_4$  (and  $c_2$ and $c_3$) are  antipodal, and } 
\endaligned}
\tag2-7
\\
\text{the mapping $(C_n\times C_n)\setminus\Cal D\ni(c,c')\mapsto c-c'\in G$ is one-to-one.}
\tag2-8
\endgather
$$
It is straightforward to verify that there exist $e_{n,i}, l_{n,i}$, $1\le i\le N_n$ such that $C_n$ satisfies  \thetag{2-6}---\thetag{2-8}.
Now let
$F_n$ be an element of  $(\Cal F_j)_{j\ge 1}$ such that $F_n\supset F_{n-1}+F_{n-1}+C_n$ and $d_{n+1}\in F_n-F_n$.
Continuing this construction process infinitely many  times we obtain a sequence $(C_n,F_{n-1})_{n\ge 1}$ satisfying \thetag{1-1}--\thetag{1-4}.
Let $T$ denote the $(C,F)$-action of $G$ associated with $(C_n,F_{n-1})_{n\ge 1}$.

{\bf Claim 1.}  $T\times T$ is ergodic.

Take $m>0$ and $v_1,v_2,w_1,w_2\in F_m$.
There is $n>m$ such that $d_n=w_2-w_1+v_1-v_2$.
Let
$$
A:=\bigsqcup_{i,j=1}^{N_n}[v_1+D+l_{n,i}]_n\times[v_2+D-e_{n,j}]_n\subset [v_1]_m\times[v_2]_m,
$$
where $D:=C_{m+1}+\cdots+ C_{n-1}$.
Define a map $\tau:A\to [w_1]_m\times[w_2]_m$ by setting
$$
\tau(x,y)=(T_{w_1-v_1+e_{n,j}-l_{n,i}}x,T_{w_1-v_1+e_{n,j}-l_{n,i}}y)
$$
if $x\in [v_1+D+l_{n,i}]_n$ and $y\in [v_2+D-e_{n,j}]_n.
$
Indeed, since 
$$
\align
T_{w_1-v_1+e_{n,j}-l_{n,i}}[v_1+D+l_{n,i}]_n&=[w_1+D+e_{n,j}]_n\text{ and}\\
  T_{w_1-v_1+e_{n,j}-l_{n,i}}[v_2+D-e_{n,j}]_n&=[w_2+D+(-l_{n,i}-d_n])_n
  \endalign
  $$ 
  for all $i,j=1,\dots N_n$, it follows that $\tau$ is  a bijection of $A$ onto $\tau (A)\subset  [w_1]_m\times[w_2]_m$.
Of course, $\tau\in[[T\times T]]$.
It is easy to compute that 
$$
(\mu\times\mu)(A)=(\mu\times\mu)([v_1]_m\times[v_2]_m)/16.
$$
By Lemma~1.1(ii), the action $T\times T$ is ergodic.

\comment
\footnote{Lemma~1.1(ii) is applied to the $(C,F)$-action $(T_g\times T_{g'})_{(g,g')\in G\times G}$ associated with the sequence $(C_{n}\times C_n,F_{n-1}\times F_{n-1})_{n\ge 1}$ and the subgroup $H=\{(g,g)\mid g\in G\}$ of $G\times G$.}

\endcomment

{\bf Claim 2}. $T\times T^{-1}$ is not ergodic.
Fix $f_1\in F_1\setminus\{0\}$.
We let
$$
Z:=\{(x,\widetilde x)\in [0]_1\times[0]_1\mid x_j\text{ and }\widetilde x_j\text{ are not antipodal for each } j>0\},
$$
where $x=(0,x_2,x_3,\dots),\widetilde x=(0,\widetilde x_2,\widetilde x_3,\dots)\in F_1\times C_2\times \cdots$.
It is easy to compute that
$$
(\mu\times\mu)(Z)=\bigg(1-4\sum_{j>1}\frac 1{N_j} \bigg)(\mu\times\mu)([0]_1\times[0]_1).
$$
Hence $(\mu\times\mu)(Z)>0$.
We show that $\bigcup_{g\in G}(T_g\times T_{-g})Z\cap([0]_1\times[f_1]_1)=\emptyset$.
If not, there is $(x,\widetilde x)\in Z$ and $g\in G$ such that $T_gx\in[0]_1$ and $T_{-g}\widetilde x\in [f_1]_1$.
Since 
 $T_gx=(0,x_2',x_3',\dots)\in F_1\times C_2\times C_3\times\cdots$ and
$T_{-g}\widetilde x=(f_1,\widetilde x_2',\dots)\in  F_1\times C_2\times\cdots$, there are integers $M_1$ and $M_2$ such that
$$
\aligned
g&=(x_2'-x_2)+\cdots+ (x_{M_1}'-x_{M_1})\quad \text{and}\\
-g&=f_1+(\widetilde x_2'-\widetilde x_2)+\cdots+ (\widetilde x_{M_2}'-
 \widetilde x_{M_2})
\endaligned
\tag 2-9
$$
with $x_{M_1}\ne x'_{M_1}$ and  $\widetilde x_{M_2}\ne \widetilde x'_{M_2}$.
It follows from \thetag{2-6}
 that  $M_1=M_2$.
 Since $x_{M_1}$ and $\widetilde x_{M_1}$ are not antipodal, we deduce from \thetag{2-7}  and \thetag{2-9} that
 $x_{M_1}-x'_{M_1}=\widetilde x_{M_1}'-\widetilde x_{M_1}$.
Hence~\thetag{2-9}
yields that
$$
\aligned
h&=(x_2'-x_2)+\cdots+ (x_{M_1-1}-x_{M_1-1}')\quad \text{and}\\
-h&=
f_1+(\widetilde x_2'-\widetilde x_2)+\cdots+ (\widetilde x_{M_1-1}'-
 \widetilde x_{M_1-1})
\endaligned
$$
where $h:=g+x_{M_1}-x'_{M_1}$.
Continuing this way several times, we find $ L\in\{2,3,\dots, M_1\}$ and $f\in G$ such that
$$
\aligned
f&=(x_2'-x_2)+\cdots+ (x_{L}'-x_L)\quad \text{and}\\
-f& =f_1+(\widetilde x_2'-\widetilde x_2)+\cdots+ (\widetilde x_{L}'-\widetilde x_{L})
\endaligned
$$
with
$x_{L}-x'_{L}\ne\widetilde x_{L}'-\widetilde x_{L}$ and $x_{L}-x'_{L}\ne 0$
and $\widetilde x_{L}'-\widetilde x_{L}\ne 0$.
If such an $L$ does not exists we then  obtain  that $f=0$ and hence $f_1=0$, a contradiction.
However then
it follows from~\thetag{2-7} that  $c_L$ and $\widetilde c_L$  are antipodal, a contradiction again.

 {\bf Claim 3.} $T\times T\times T$ is not conservative.
 Let
 $$
 W:=\{(x,y,z)\in[0]_0\times[0]_0\times[0]_0\mid x_j\ne y_j,y_j\ne z_j, x_j\ne z_j\text{ for each $j>0$}\},
 $$
 where $x_j,y_j$ and $z_j$ are the $j$-th coordinate of $x,y$ and $z$ viewed as infinite sequences from $\{0\}\times C_1\times C_2\times\cdots$.
 Then
 $$
 (\mu\times\mu\times\mu)(W)>1-\frac 3 4\sum_{j>0}\frac 1 {N_j}>0.
 $$
 We claim that $W$ is a wandering subset for $T\times T\times T$.
 If not, there is $(x,y,z)\in W$ and $g\in G$ such that
 $$
 T_gx,T_{g}\widetilde y,
T_{g}\widetilde z\in [0]_1.\tag2-10
$$
We write the expansions  $x=(0,x_1,x_2,\dots)$,
$y=(0,\widetilde y_1,\widetilde y_2,\dots)$,
$z=(0,\widetilde z_1,\widetilde z_2,\dots)$,
 $T_gx=(0,x_1',x_2',\dots)$,
$T_{g}y=(0,\widetilde y_1',\widetilde y_2',\dots)$ and
$T_{g}z=(0,\widetilde z_1',\widetilde z_2',\dots)$ as infinite sequences from $ \{0\}\times C_1\times C_2\times\cdots$.
Then \thetag{2-10}  and \thetag{2-6} yield that there is integer $M$  such that
$$
\aligned
g&=(x_1'-x_1)+\cdots+ (x_{M}'-x_{M}),\\
g&= (y_1'-y_1)+\cdots+ (y_{M}'-
 y_{M})
 \quad \text{and}\\
 g&=(z_1'-z_1)+\cdots+ (z_{M}'-
 z_{M})
\endaligned
\tag2-11
$$
with $x_{M}'-x_{M}\ne 0$, $y_{M}'-
 y_{M}\ne 0$ and $z_{M}'-
 z_{M}\ne 0$.
If 
 $x_{M}'-x_{M}= y_{M}'-
 y_{M}$ then $x_M=y_M$ by \thetag{2-8} and hence $(x,y,z)\not\in W$.
 Therefore $x_{M}'-x_{M}\ne  y_{M}'-
 y_{M}$.
 In a similar way,
 $y_{M}'-
 y_{M}\ne z_{M}'-z_{M}$.
 However then \thetag{2-11} and \thetag{2-7} yield that
  $x_{M_1}$ and $y'_{M_1}$ are antipodal
and  $z_{M_1}$ and $y'_{M_1}$ are antipodal.
 This is only possible if  $x_{M_1}=z_{M_1}$ and hence  $(x,y,z)\not\in W$,
 a contradiction.
 \comment

 In the case $G=\Bbb Z$, we argue in a similar way but set
  $$
C_n:=\{a-e_i,a+e_i,a+l_i, a-l_i-d_n\mid i=1,\dots N_n\},
$$
where $e_i,  l_i$, $1\le i
\le N_n$, and $a$ are  are  positive integers   chosen in such a way  that
\thetag{2-1}  and \thetag{2-2} are satisfied and $\min C_n=\{0\}$.
We also set
$$
F_n:=\{0,1,\dots, \max(F_{n-1}+C_n)\}.
$$
\endcomment

The fact that $T\times T^{-1}$ is conservative follows from Theorem~0.4.
 \qed
 \enddemo

\comment

\demo{Proof of Theorem~{0.4}}
This  is a slight modification of the proof of Theorem~0.3.
That is why our exposition is rather sketchy. 
Suppose that we have already determined  the sequence $(C_j,F_{j})_{j=1}^{n-1}$.
Select an element $d_n\in F_{n-1}-F_{n-1} $ and  set
$$
C_n:=\{e_{n,i},l_{n,i}, L_n+e_{n,i}, L_n+ l_{n,i}-d_n\mid i=1,\dots N_n\},
$$
for some elements $L_n, e_{n,i},  l_{n,i}$ of $G$, $1\le i
\le N_n$ such that \thetag{2-5} holds.
We call
 $e_{n,i}$ and $L_n+l_{n,i}$ as well as  $l_{n,i}$ and $L_n+l_{n,i}-d_n$ {\it antipodal}, $1\le i\le N_n$.
It is straightforward to verify that there exist $L_n,e_{n,i}, l_{n,i}$, $1\le i\le N_n$ such that $C_n$ satisfies both \thetag{2-5} and \thetag{2-6}.
Now let
$F_n$ be an element of the sequence $(\Cal F_j)_{j\ge 1}$ such that $F_n\supset F_{n-1}+C_n$.
Continuing this construction process infinitely many  times we obtain a sequence $(C_n,F_{n-1})_{n\ge 1}$ satisfying \thetag{1-1}--\thetag{1-4}.
We will assume in addition that $(d_n)_{n=1}^\infty$ is a sequence in which each element of $G$ occurs infinitely many times
and that  $\sum_{n>0}\frac 1{N_n}<\frac 14$.

Let $T$ denote the $(C,F)$-action of $G$ associated with $(C_n,F_{n-1})_{n\ge 1}$.

{\bf Claim 1.}  $T\times T^{-1}$ is ergodic.

Take $m>0$ and $v,w\in F_m$.
There is $n>m$ such that $d_n=w+v$.
Let
$$
A:=\bigsqcup_{i,j=1}^{N_n}[D+l_{n,i}]_n\times[D+L_n+e_{n,j}]_n\subset [0]_m\times[0]_m,
$$
where $D:=C_{m+1}+\cdots+ C_{n-1}$.
Define a map $\tau:A\to [v]_m\times[w]_m$ by setting
$$
\tau(x,y)=(T_{v+e_{n,j}-l_{n,i}}x,T_{-v-e_{n,j}+l_{n,i}}y)\quad\text{if}\quad x\in [D+l_{n,i}]_n\quad \text{and }y\in [D-e_{n,j}]_n.
$$
Since 
$$
\align
T_{v-l_{n,i}+e_{n,j}}[D+l_{n,i}]_n&=[v+D+e_{n,j}]_n\text{ and}\\
  T_{-v+l_{n,i}-e_{n,j}}[D+L_n+e_{n,j}]_n&=[w+D+(L_n+l_{n,i}-d_n)]_n
  \endalign
  $$ 
  for all $i,j=1,\dots N_n$, it follows that $\tau$ is  a bijection of $A$ onto $\tau (A)\subset  [v]_m\times[w]_m$.
Extending $\tau$ in a natural way to the entire space, we obtain $\tau\in[T\times T^{-1}]$ and $(\mu\times\mu)(A)=(\mu\times\mu)([0]_m\times[0]_m)/16$.
By Lemma~1.1(ii), the action $T\times T^{-1}$ is ergodic.

{\bf Claim 2}. $T\times T$ is not ergodic.

Fix $d_1\in F_1\setminus\{0\}$.
Define a subset $Z\subset[0]_1\times[0]_1$  in the same way as in Claim~2 of the proof of Theorem~0.3.
It is of positive measure.
We show that the $(T\times T)$-orbit of $Z$ does not intersect the cylinder $[0]_1\times[d_1]_1$
If not, there is $(x,\widetilde x)\in Z$ and $g\in G$ such that $T_gx\in[0]_1$ and $T_{g}\widetilde x\in [d_1]_1$.
Since $x=(0,c_2,\dots)\in  F_1\times C_2\times\cdots$,
$\widetilde x=(0,\widetilde c_2,\dots)\in  F_1\times C_2\times\cdots$,
 $T_gx=(0,c_2',\dots)\in  F_1\times C_2\times\cdots$ and
$T_{g}\widetilde x=(d_1,\widetilde c_2',\dots)\in  F_1\times C_2\times\cdots$,
it follows from \thetag{2-5}
that  there is  $M_1$  such that
$$
\aligned
g&=(c_2'-c_2)+\cdots+ (c_{M_1}'-c_{M_1})\quad \text{and}\\
g&=d_1+(\widetilde c_2'-\widetilde c_2)+\cdots+ (\widetilde c_{M_1}'-
 \widetilde c_{M_1})
\endaligned
\tag 2-8
$$
with $c_{M_1}\ne c'_{M_1}$ and  $\widetilde c_{M_1}\ne \widetilde c'_{M_1}$.
 Since $c_{M_1}$ and $\widetilde c_{M_1}$ are not antipodal, we deduce from \thetag{2-6} that
 $c'_{M_1}-c_{M_1}=\widetilde c_{M_1}'-\widetilde c_{M_1}$.
Now~\thetag{2-8}
yields that
$$
\aligned
h&=(c_1'-c_1)+\cdots+ (c_{M_1-1}-c_{M_1-1}')\quad \text{and}\\
h&=
d_1+(\widetilde c_1'-\widetilde c_1)+\cdots+ (\widetilde c_{M_2-1}'-
 \widetilde c_{M_2-1})
\endaligned
$$
where $h=g+c_{M_1}-c'_{M_1}$.
Continuing this way several times, we find $L\le M_1$ and $f\in G$ such that
$$
\aligned
f&=(c_1'-c_1)+\cdots+ (c_{L}'-c_L)\quad \text{and}\\
f& =d_1+(\widetilde c_1'-\widetilde c_1)+\cdots+ (\widetilde c_{L}'-\widetilde c_{L})
\endaligned
$$
with
$c_{L}-c'_{L}\ne\widetilde c_{L}'-\widetilde c_{L}$ and $c_{L}-c'_{L}\ne 0$
and $\widetilde c_{L}'-\widetilde c_{L}\ne 0$.
Of course, $L\ge 1$ because otherwise $f=0$ and hence $d_1=0$.
However then
it follows from~\thetag{2-6} that  $c_L$ and $\widetilde c_L$  are antipodal, a contradiction.

 {\bf Claim 3.} $T\times T$ is conservative.
 
 Let $n$ be such that $d_n=0$.
 We set $A:=\bigsqcup_{i,j=1}^{N_n}[e_{n,i}]_n\times[l_{n,j}]\subset [0]_{n-1}\times[0]_{n-1}$.
 Then $(\mu\times\mu)(A)=(\mu\times\mu)([0]_{n-1}\times[0]_{n-1})/16$ and
 $
 (T_{L_n}\times T_{L_n})A\subset [0]_{n-1}\times[0]_{n-1}.
 $
 By Lemma~1.1(i), $T\times T$ is conservative.

 {\bf Claim 4.} $T\times T\times T$ is not conservative.
 This is proved almost literally as Claim~3 from the proof of Theorem~0.3.
  \qed
 \enddemo

\endcomment

 \demo{Proof of Theorem~0.4}
 For each $n>0$, we
 let
 $$
 A:=\bigsqcup_{c\ne c'\in C_{n+1}}[c]_{n+1}\times[c']_{n+1}\subset [0]_n\times[0]_n
 $$
and  define a map $\tau: A\to [0]_n\times[0]_n$ by setting
 $$
 \tau(x,y)=(T_{c'-c}x, T_{c-c'}y)\quad\text{if $x\in[c]_{n+1},y\in [c']_{n+1}$}.
 $$
Then $\tau([c]_{n+1}\times[c']_{n+1})=[c']_{n+1}\times[c]_{n+1}$, $\tau\in[[T\times T^{-1}]]$ and
 $$
 (\mu\times\mu)(A)=\bigg(1-\frac 1{\# C_{n+1}}\bigg)(\mu\times\mu)([0]_n\times[0]_n)\ge \frac 12(\mu\times\mu)([0]_n\times[0]_n).
 $$
 By Lemma~1.1(i), $T\times T^{-1}$ is conservative.
 \qed
 \enddemo

\comment
Given $n>0$,  $l,i\in\{1,\dots,k\}$ and $\bold c\ne \bold c'\in C_{n+1}^{l,i}$, there is an element $g(\bold c,\bold c')\in G\setminus\{0\}$
such that $g(\bold c,\bold c')+\bold c=\bold c'$.
We now set
$$
A:=\bigsqcup_{i=1}^k\bigsqcup_{\bold c\ne \bold c'\in C_{n+1}^{l,i}}[(0,l)*\bold c]_{n+1}\times[(0,l)*\bold c']_{n+1}\subset[(0,l)]_n\times [(0,l)]_n
$$
and define a map $\tau: A\to [(0,l)]_n\times[(0,l)]_n$ by setting
 $$
 \tau(x,y)=(T_{g(\bold c,\bold c')}x, T_{-g(\bold c,\bold c')}y)\quad\text{if $x\in[(0,l)*\bold c]_{n+1},y\in [(0,l)*\bold c']_{n+1}$}.
 $$
Of course, $\tau$ extends naturally to a transformation from the full group $[T\times T^{-1}]$ and
$$
 (\mu\times\mu)(A)=\sum_{i=1}^k\bigg(1-\frac 1{r_{n+1}^{l,i}}\bigg)r_{n+1}^{l,i}r_{n+1}^{l,i}\lambda_{n+1}^i\lambda_{n+1}^i
 $$
 
 \endcomment

\head 3. On ``symmetry'' of Markov shifts
\endhead

In this section we consider only the case where $G=\Bbb Z$.
We first recall some basic definitions and properties of  infinite measure preserving Markov shifts.

Let $S$ be an infinite countable set and let $P=(P_{a,b})_{a,b\in S}$ be a stochastic matrix over $S$.
Suppose that there is a strictly positive vector $\lambda=(\lambda_s)_{s\in S}$ which is  a left eigenvector for $P$ with eigenvalue $1$, i.e. $\lambda P=\lambda$.
Moreover, assume that $\sum_{s\in S}\lambda_s=\infty$.
Consider the infinite product space $X:=S^\Bbb Z$ and endow $X$ with the infinite product Borel structure.
Let $T$ denote the left shift on $X$.
Given $s_0,\dots,s_{k}\in S$ and $l\in\Bbb Z$, we denote by $[s_0,\dots,s_k]_l^{l+k}$ the cylinder
$\{x=(x_j)_{j\in \Bbb Z}\mid x_l=s_0,\dots, x_{l+k}=s_k\}$.
Define a  measure $\mu_{P,\lambda}$ on  $X$ by setting $\mu_{P,\lambda}([s_0,\dots,s_k]_l^{l+k})=\lambda_{s_0}P_{s_0,s_1}\cdots P_{s_{k-1},s_k}$ for each cylinder $[s_0,\dots,s_k]_l^{l+k}$.
Then $\mu_{P,\lambda}$ extends uniquely to the Borel $\sigma$-algebra on $X$ as a $\sigma$-finite
infinite   measure which is invariant under $T$.
The dynamical system $(X,\mu_{P,\lambda},T)$ is called an  {\it infinite Markov shift}.

\proclaim{Lemma 3.1 (\cite{Aa}, \cite{KaPa})}
$(X,\mu_{P,\lambda},T)$ is conservative and ergodic if and only if the following two conditions are satisfied:
\roster
\item"$(i)$"  $P$ is irreducible, i.e. for each $a,b\in S$, there is $n>0$ such that $P^{(n)}_{a,b}>0$ and
\item"$(ii)$"
$P$ is recurrent, i.e. $\sum_{n>0}P^{(n)}_{a,a}=\infty$ for some (and hence for each in view of $(i)$) $a\in S$.
\endroster
If $(ii)$ does not hold then $(X,\mu_{P,\lambda},T)$ is not conservative.
\endproclaim

Here $P^{(n)}$ means the usual matrix power $P\cdots P\,(n\text{ times})$.

Let $\sigma:X\to X$ denote the flip, i.e. $(\sigma x)_n:=x_{-n}$ for $x\in X$ and $n\in\Bbb Z$.
Denote by $\Lambda=(\Lambda_{a,b})_{a,b\in S}$ a matrix over $S$ such that $\Lambda_{a,b}=\lambda_a$ if $a=b$ and 
$\Lambda_{a,b}=0$ if $a\ne b$.
It is straightforward to verify  that $\sigma T\sigma^{-1}=T^{-1}$, $\Lambda^{-1} P^*\Lambda$ is a stochastic matrix and $\mu_{P,\lambda}\circ\sigma=\mu_{\Lambda^{-1} P^*\Lambda,\lambda}$.
Given two infinite Markov shifts which are defined on the spaces $(S^{\Bbb Z},\mu_{P,\lambda})$ and $(S_1^{\Bbb Z},\mu_{P_1,\lambda_1})$, their Cartesian product is an infinite  Markov shift  defined on the space $((S\times S_1)^\Bbb Z,\mu_{P\otimes P_1,\lambda\times\lambda_1})$, where the matrix $P\otimes P_1$ is defined over $S\times S_1$ by  $(P\otimes P_1)_{(a,a_1),(b,b_1)}:=P_{a,b}P_{b,b_1}$.

\proclaim{Corollary 3.2} Let $(X,\mu_{P,\lambda},T)$ be an infinite  Markov shift and let $0\le m\le k$.
Then the transformation 
 $\underbrace{T\times\cdots\times T}_{m \text{ times}}\times \underbrace{T^{-1}\times\cdots\times T^{-1}}_{k-m \text{ times}}$ is  conservative and ergodic if and only if $T\times\cdots\times T( k\text{ times})$ is  conservative and ergodic.
\endproclaim

\demo{Proof} Fix $a\in S$.
Then
$$(P^{\otimes m}\otimes (\Lambda^{-1} P^*\Lambda)^{\otimes (k-m)})^{(n)}_{(a,\dots,a)}=(P^{(n)}_{a,a})^m
(P^{(n)}_{a,a})^{k-m}=(P^{(n)}_{a,a})^k=(P^{\otimes k})_{(a,\dots, a)}^{(n)}.
$$
Hence by Lemma~3.1(ii), the stochastic matrix $P^{\otimes m}\otimes (\Lambda^{-1} P^*\Lambda)^{\otimes (k-m)}$
is recurrent if and only if  the stochastic matrix $P^{\otimes k}$ is recurrent.
In a similar way one can check that $P^{\otimes m}\otimes (\Lambda^{-1} P^*\Lambda)^{\otimes (k-m)}$ is irreducible if and only if so is  $P^{\otimes k}$.
\qed
\enddemo

The following assertion follows from Lemma~3.1 and Corollary~3.2.

\proclaim{Corollary 3.3} Let $T$ an ergodic  conservative infinite Markov shift of ergodic index one.
Then $T\times T^{-1}$ is not conservative.
\endproclaim

We note that Corollary~3.3 was proved in \cite{Cl--Va} under an extra assumption
that $P=\Lambda^{-1} P^*\Lambda$.

\head 3. Open problems and remarks
\endhead

\roster
\item
Given $p\ge k\ge 1$, is there a mixing rank-one infinite measure preserving transformation of ergodic index $k$ such that
$T\times\cdots \times T(l\text{ times})$ is conservative if and only if $l\le p$? 
Theorem~0.2 provides an affirmative answer to this question if $p=k+1$.
\item
Is there a rank-one infinite measure preserving transformation $T$ such that   $T\times T^{-1}$ is ergodic but $T\times T$ is not?
\item 
Is there a rank-one  infinite measure preserving transformation $T$ such that $T\times T\times T$ is ergodic but $T\times T^{-1}$ is not?
\item
We note that Theorem~0.4 extends naturally to the  ergodic infinite measure preserving actions  of {\it finite funny rank} (see \cite{Da2} for the definition). 
\item
It would be interesting to investigate which indexes of ergodicity and conservativeness are realizable on the infinite measure preserving transformations which are Maharam extensions of type $III_1$ ergodic nonsingular transformations (see \cite{DaSi} for the definitions).
\endroster

\enddocument